\documentclass{amsart}
\usepackage{amssymb, latexsym, amsmath}%, makeidx}
\usepackage[all]{xy}
\usepackage{enumerate}
\usepackage[usenames, dvipsnames]{color}
\usepackage{mathrsfs} %nice script font
\usepackage[OT2,T1]{fontenc}
\DeclareSymbolFont{cyrletters}{OT2}{wncyr}{m}{n}
\DeclareMathSymbol{\Sha}{\mathalpha}{cyrletters}{"58}

\theoremstyle{plain}
\newtheorem{Theorem}[equation]{Theorem}
\newtheorem{Lemma}[equation]{Lemma}
\newtheorem{Corollary}[equation]{Corollary}

\newtheorem*{Conjecture*}{Conjecture}

\theoremstyle{definition}

\newtheorem{Paragraph}[equation]{}

\theoremstyle{remark}
\newtheorem{Remark}[equation]{Remark}

\numberwithin{equation}{section}

\renewcommand{\b}{{\text{\rm b}}}
\newcommand{\h}{{\text{\rm h}}}
\renewcommand{\k}{{\kappa}}

\newcommand{\D}{{\mathscr D}}

\newcommand{\G}{{\text{\rm G}}}

\renewcommand{\H}{{\text{\rm H}}}

\renewcommand{\O}{{\text{\rm O}}}

\newcommand{\Q}{{\mathbb Q}}

\renewcommand{\S}{{\mathcal S}}

\newcommand{\Z}{{\mathbb Z}}

\newcommand{\br}[1]{{\left<{#1}\right>}}

\newcommand{\car}{{\text{\rm char}}}

\newcommand{\cs}{{\text{\rm cs}}}
\renewcommand{\det}{{\text{\rm det}}}
\newcommand{\df}{{\,\overset{\text{\rm df}}{=}\,}}

\renewcommand{\div}{{\text{\rm div}}}

\newcommand{\im}{{\text{\rm im}}}
\newcommand{\ind}{{\text{\rm ind}}}
\renewcommand{\inf}{{\text{\rm inf}\,}}

\newcommand{\isom}{{\;\simeq\;}}
\renewcommand{\ker}{{\text{\rm ker}}}
\newcommand{\lcm}{{\text{\rm lcm}}}

\newcommand{\ns}{{\text{\rm ns}}}

\newcommand{\red}{{\text{\rm red}}}

\newcommand{\res}{{\text{\rm res}}}
\newcommand{\rk}{{\text{\rm rk}}}

\newcommand{\Br}{{\text{\rm Br}}}

\newcommand{\Card}{{\text{\rm Card}}}

\newcommand{\Div}{{\text{\rm Div}\,}}

\newcommand{\Frac}{{\text{\rm Frac}\,}}

\newcommand{\Supp}{{\text{\rm Supp}}}
\newcommand{\Spec}{{\text{\rm Spec}\,}}

\begin{document}
%%   Here starts the topmatter.
%% Text must not happen before \maketitle command!

\title[Cyclicity over $p$-Adic Curves]
{Cyclicity in Prime Degree over a $p$-Adic Curve}
%%   NEW! Title now has optional argument like other sectioning commands.
%% It will appear in running heads on even pages exept for the first one.
%%   The second mandatory argument produces the full title of your paper.
%% It will appear in running head in case you omit the optional one.
%% This title consists of two lines, add extra \\ to get more lines.
%% We think it is a good idea to avoid math formulas in the title.
%% If any they will be boldface.

\author{Eric Brussel}
%, Kelly McKinnie, and Eduardo Tengan}
\address 
{Department of Mathematics\\
California Polytechnic State University\\
San Luis Obispo, CA 93407\\ USA}
\email{ebrussel@calpoly.edu}

%11G20 is curves over local and finite fields
%11R58 is arithmetic theory of function fields
%14G40 is arithmetic varieties and schemes; Arakolov theory and heights.
%14E22 is ramification problems
%16K50 is Brauer groups
%12G05 is Galois cohomology
%16K20 is fd associative algebras
%16K50 is Brauer groups.
%16S35 is twisted and skew group rings, crossed products.

%%   This is the example of two-authors article. Separate them by \and.
%%   NEW! Function of \address and \email is similar to that in amsppt.
%% But the usage is different: here you place them INSIDE the argument
%% of \author. \address and \email produce new line automatically, so
%% don't put \\ before them.
%%   Note that if you have thanks, they should go exactly after the
%% authors name preceding the \address. If your paper has many authors,
%% please try to make the same number of lines of the address for all
%% authors.

%\dedicatory{Include a dedication?}
%%   NEW! Dedication is a best way to impress your gratitudes and sympathies.

%\received?} \accepted{September 22, 2010}
%%%   It's NEW! but you'll never deal with it, it's a headache of
%%% your publisher.

\begin{abstract}
We reprove two results of Saltman, \cite[Theorem 5.1, Corollary 5.2]{Sa07}:
If $F$ is the function field of a smooth $p$-adic curve and
$D$ is an $F$-division algebra of prime degree $\ell\neq p$ then $D$ is $\Z/\ell$-cyclic,
and that if $D$ is an $F$-division algebra of prime period $\ell\neq p$ then $D$ has index
$\ell$ if and only if its ramification locus on a suitable 2-dimensional model for $F$ has no ``hot points''.
\end{abstract}
%%   We insist on abstract in your paper!

\hfill February 4, 2014
\maketitle

%\section*{Introduction}
%%   Use unstarred sectioning commands to obtain automatic numeration.
%% And don't type dots at the end of a heads.
%%   NEW! The number of accessible sectioning commands was highly restricted.
%% Here is the whole list: \section \subsection \subsubsection. This keeps
%% you from making the structure of text too complex.
\section*{Introduction}

One of the most important open problems in the area of finite-dimensional division algebras 
is to determine if every division algebra of prime degree over a field is cyclic
(see \cite[Section 1]{ABGV}).
Attempts to solve it often involve analyses over fields $F$ for which there is 
a reasonable theory of arithmetic.
In the seminal result \cite[Theorem 5.1]{Sa07}
Saltman solved the problem for $F$ the function field of a smooth $p$-adic curve (e.g., $F=\Q_p(t)$),
proving that all $F$-division algebras of prime degree $\ell\neq p$ are $\Z/\ell$-cyclic.
We reprove this result as a corollary of a slight generalization, which is that if 
$F$ is the function field of a smooth curve over a complete discretely valued field $K=(K,v)$,
$\car(\k(v))=p\geq 0$,
and $D$ is an $F$-division algebra of prime degree $\ell\neq p$, then there exists a $\Z/\ell$-cyclic
field extension $L/F$ such that $\alpha_L$ is unramified.
In \cite[Theorem 7.13]{Sa08} Saltman proved this generalization for an arbitrary regular surface under
the assumption that $F$ contain the $\ell$-th roots of unity.
We also reprove \cite[Corollary 5.2]{Sa07},
which states that if $F$ is the function field of a $p$-adic curve
and $D$ is an $F$-division algebra of prime period $\ell\neq p$ then $D$ has index
$\ell$ if and only if its ramification locus on a suitable 2-dimensional model for $F$ has no ``hot points''.

We use the machinery and methods of \cite{BMT13}, which we view as a kind of extension of Grothendieck's proper
base change theorem in the following sense.  
Let $F$ be the function field of a smooth curve over a complete discretely valued field $K=(K,v)$.
Then $F$ is the function field of a regular relative curve $X/\O_v$, 
and the reduced scheme $C$ underlying the closed
fiber $X\otimes_K\k(v)$ is a projective curve over $\k(v)$.
Let $n$ be a number prime to $\car(\k(v))$, and let $\H^q(-)$ denote an \'etale cohomology group with
$n$-torsion coefficients.
By Grothendieck's theorem we have isomorphisms $\H^q(X)\isom\H^q(C)$ in all degrees $q$, whereby elements
of $\H^q(C)$ ``lift'' to $\H^q(X)$.
In \cite{BMT13} we showed how to extend this lifting to a subgroup of $\H^q(\k(C))$, resulting in 
constructions of elements of $\H^q(F)$ with controlled ramification behavior.
By manipulating the model $X$ (using blow ups) 
we can smooth out the ramification divisor of a given element of $\H^q(\k(C))$,
to a point where it is within reach of our lift.
In the current paper we show that if $\alpha\in\H^2(F,\mu_\ell)$ is of prime degree $\ell\neq p$
then there exists a model $X$ over which $\alpha$'s ramification divisor
is subject to splitting by a $\Z/\ell$-cyclic extension of $F$ that is lifted
from a cyclic extension of $\k(C)$ 
constructed using Saltman's generalized Grunwald-Wang theorem \cite[5.10]{Sa82}.

Saltman took a more overtly valuation-theoretic approach in \cite{Sa07}, manipulating the model $X$
until he could define an element $f\in F$ with divisor 
$\div(f)$ approximating the division algebra's ramification divisor, and then descending
the cyclic extension $F(\mu_\ell)(f^{1/\ell})$ to $F$.  
Instead of trying to copy the ramification divisor, our strategy is essentially
to manipulate $X$ so that we can glue together the residues, which are 
(tamely ramified) cyclic covers of the ramification divisor's prime factors.

\section{Notation and Background}

\Paragraph{\bf General Conventions.}
Let $S$ be an excellent scheme, $n$ a number that is invertible on $S$,
and $\Lambda=(\Z/n)(i)$ the \'etale sheaf $\Z/n$ twisted by an integer $i$.
We write $\H^q(S,\Lambda)$ for the \'etale cohomology group,
and if $\Lambda$ is arbitrary and fixed (or doesn't matter) 
we write $\H^q(S)$ instead of $\H^q(S,\Lambda)$,
and $\H^q(S,r)$ in place of $\H^q(S,\Lambda(r))$.
If $S=\Spec A$ for a ring $A$, then we 
write $\H^q(A)$ and $\H^q(A,r)$.
If $v$ is a valuation on a field $F$, we write $\k(v)$ for
the residue field of the valuation ring $\O_v$, and $F_v$ for the completion of $F$
at $v$.  If $v$ arises from a prime divisor $D$ on $S$, we write $v=v_D$, $\k(D)$,
and $F_D$.
If $T$ is an integral closed subscheme of $S$ we write $\k(T)$ for its function field. 
If $T\to S$ is a morphism of schemes, then
the restriction $\res_{S|T}:\H^q(S)\to\H^q(T)$ is defined,
and we write $\beta|_T$ or $\beta_T$ for $\res_{S|T}(\beta)$.

\Paragraph{\bf Basic Setup.}\label{setup}
Let $R$ be a complete discrete valuation ring with fraction
field $K$ and residue field $k$, $n$ a number invertible in $R$, 
$F$ the function field of a smooth projective curve over $K$,
and $X/R$ a regular (projective, flat) relative curve with function field $F$.
Thus $X$ is 2-dimensional, all of its closed points $z$ have codimension 2 (\cite[8.3.4]{Liu}), 
and the corresponding local rings $\O_{X,z}$ are factorial.
The closed fiber $X_0=X\otimes_R k$ is a connected projective curve over $k$; 
write $C=X_{0,\red}$ for its reduced subscheme, 
$C_1,\dots,C_m$ for the irreducible components of $C$, and $\k(C)=\prod_i\k(C_i)$.
We assume throughout that each irreducible component $C_i$ of $C$ is regular, and that all
singular points of $C$ have multiplicity two, a situation that can always be achieved by blowing up.
We let $\S$ denote the set of singular points of $C$.
If $z$ is a closed point of $X$ then $z$ lies on $C$,
and if $z\in C_i$ we write $K_{i,z}=\Frac(\O_{C_i,z}^\h)$, 
where the superscript ``$\h$'' denotes {\it henselization}. 
If $z\in\S$ is on $C_i\cap C_j$ we let $K_z=\Frac(\O_{C,z}^\h)=K_{i,z}\times K_{j,z}$.

Since exactly two irreducible components of $C$ meet at any $z\in\S$
the {\it dual graph} $\G_C$ is defined, and consists of a vertex for each irreducible
component of $C$ and an edge for each singular point, 
such that an edge and a vertex are incident when the corresponding singular point
lies on the corresponding irreducible component (\cite[2.23]{Sai85}, see also \cite[10.1.48]{Liu}). 
The (first) {\it Betti number} $\b_1(C)=\rk(\H_1(\G_C,\Z))$ is the sum $N+E-V$,  
where $V,E$ and $N$ are the numbers of vertices, edges, and connected components of 
$\G_C$.
Since $C$ is connected we have $N=1$, and $\b_1(C)=1-m+s$, where $s=\Card(\S)$,
and this number computes the number of chordless cycles of $\G_C$.
In particular it is zero if $\G_C$ is a tree.

%Let $\chi(C)$ denote the {\it Euler characteristic} of the dual graph $\Gamma_C$ of $C$.
%We have $\chi(C)=1-\b_1(C)=m-s$.  
%See Amer. Math. Monthly 78 no.4, 1971, How is a Graph's Betti Number Related to Its Genus?
%by R.A. Duke.  $\chi(C)$ equals V-E, which is also 1-(number of cycles).

We call an effective (Cartier) divisor $D\in\Div X$ {\it regular} if its underlying closed subscheme 
$\Supp(D)$ is regular,
and {\it horizontal} if $\Supp(D)$ is finite over $\Spec R$.
Assume $C$ is as above, with regular irreducible components and singular points of
multiplicity two.
By \cite[1.12]{BMT13} it is possible to choose for each closed point $z\in C_{(0)}$ a
{\it distinguished prime divisor} $D_z$ on $X$,
which is a regular horizontal prime divisor that passes through $z$
and is transverse to each irreducible component of $C$ passing through $z$.
Let $\D$ denote a set of distinguished prime divisors, one for each $z\in C_{(0)}$.
Let $\D_*\subset\D$ denote those distinguished prime divisors $D_z$ with $z\in\S$,
and let $\D_\S=\D-\D_*$.  We will sometimes say a divisor is ``in $\D$'' if it is composed
of distinguished prime divisors.
Finally, if $D\in\Div X$ is any divisor we write $\D(D)$ or $\D_\S(D)$ for the support of
$D$ that is in $\D$ or $\D_\S$.

\Paragraph{\bf Residues and Ramification.}\label{residues}
All valuations will be discrete of rank one.
If $F=(F,v)$ is a discretely valued field, $F_v$ is the completion of $F$ at $v$,
$n$ is prime to $\car(\k(v))$, and $\Lambda=(\Z/n)(i)$ for some $i$,
then the {\it residue map} $\partial_v$ is defined by the diagram
\[\xymatrix{
& &\H^q(F)\ar[d]\ar@{-->}[dr]^-{\partial_v}\\
0\ar[r]&\H^q(\k(v))\ar[r]&\H^q(F_v)\ar[r]&\H^{q-1}(\k(v),-1)\ar[r]&0}\]
We call the bottom row a {\it Witt exact sequence}.
The bottom surjection is split (non canonically) by the map $\omega\mapsto(\pi_v)\cdot\omega$,
where $(\pi_v)\in\H^1(F_v,\mu_n)$ is the Kummer element
determined by a choice of uniformizer $\pi_v$ for $F_v$.
We call the resulting direct sum decomposition $\H^q(F_v)\isom\H^q(\k(v))\oplus\H^{q-1}(\k(v),-1)$
a {\it Witt decomposition}.
We say $\alpha\in\H^q(F)$ is {\it unramified at $v$}
if $\partial_v(\alpha)=0$, and {\it ramified at $v$} if $\partial_v(\alpha)\neq 0$.
If $\alpha$ is unramified at $v$ then $\alpha$ comes from the subgroup $\H^q(\O_v)\leq\H^q(F)$
(see \cite[Section 3]{CT95}),
and we say $\alpha$ is {\it defined} at $\k(v)$.
If $\alpha$ is defined at $\k(v)$ then it has a {\it value} 
\[\alpha(v)=\res_{\O_v|\k(v)}(\alpha)\;\in\H^q(\k(v))\]
Since the ring homomorphism $\O_v\to\k(v)$ factors through the complete discrete valuation
ring $\O_{F_v}\subset F_v$ we have
the alternative description $\alpha(v)=\res_{F|F_v}(\alpha)$, using the Witt sequence
to identify $\H^q(\k(v))$ with the subgroup $\H^q(\O_{F_v})\leq\H^q(F_v)$.
Note $\alpha_{F_v}=0$ if and only if
$\partial_v(\alpha)=0$ and $\alpha(v)=0$ by the Witt sequence.
If $v$ arises from a prime divisor $D$ on an integral scheme
with function field $F$ we generally substitute $D$ for $v$,
and write $\partial_D$
and $\alpha(D)$ in place of $\partial_v$ and $\alpha(v)$.

Each $f\in F^*$ defines a Kummer element $(f)\in\H^1(F,\mu_n)$,
and $\partial_v(f)=v(f)\pmod n\in\H^0(\k(v),\Z/n)$.
If $\chi\in\H^1(F,\Z/n)$ we write $(f)\cdot\chi\in\H^2(F,\mu_n)$ for the cup product.  
Then 
\begin{equation}\label{cupresidue}
\partial_v((f)\cdot\chi)=\left[v(f)\cdot\chi-(f)\cdot\partial_v(\chi)+(-1)\cdot v(f)\cdot\partial_v(\chi)\right]_{F_v}
\end{equation}
See \cite[II.7.12, p.18]{GMS} for the general cup product formula.

In the situation of \eqref{setup}
any prime divisor $D\in\Div X$ determines a valuation $v_D$ on $F$.
Thus for each $\alpha\in\H^q(F)$ we have a {\it ramification divisor}
\[D_\alpha=\sum_{\underset{\mbox{\scriptsize $\partial_D(\alpha)\neq 0$}}{D\in\Div X}}D\;\;\in\Div X\]
For a fixed $\alpha$ we may always blow up $X$ until the horizontal prime factors of $D_\alpha$ 
are all regular, avoid $\S$, and intersect each irreducible component of $C$ transversely.

In the situation of \eqref{setup} Kato defines a complex
\begin{equation}\label{kato}
\xymatrix{\H^2(F,\mu_n)\ar[r]^-{\partial_2}&\bigoplus_{D\in\Div X}\H^1(\k(D),\Z/n)\ar[r]^-{\partial_1}
&\bigoplus_{z\in C_{(0)}}\H^0(\k(z),\mu_n^*)
}\end{equation}
where $\mu_n^*=(\Z/n)(-1)$, $\partial_2=\sum_D\partial_D$,
and $\partial_1=\sum_z(\bigoplus_D\partial_{D,z})$ is composed of residue maps 
$\partial_{D,z}:\H^1(\k(D),\Z/n)\to\H^0(\k(z),\mu_n^*)$ defined in \cite{Ka86} to include the
possibility that $z$ is a singular point of $D$.
In particular if $\alpha\in\H^2(F)$ and exactly two components $D_1,D_2\subset D_\alpha$ 
run through a point $z$, then 
\begin{equation}\label{sumtozero}
\partial_z(\partial_{D_1}(\alpha))+\partial_z(\partial_{D_2}(\alpha))=0
\end{equation}

\Paragraph{\bf Hot and Cold Points.}\label{hotpt}
Assume the situation of \eqref{setup}.
Suppose $n=\ell\neq p$ is a prime, $\alpha\in\H^2(F,\mu_n)$, and $D_\alpha$ has regular irreducible
components, at most two of which meet at any given point.
Following Saltman in \cite{Sa07}
we will say a singular point $z\in D_\alpha$ at the intersection of $D_1,D_2\subset D_\alpha$ 
is a {\it hot point} if each $\theta_{D_i}=\partial_{D_i}(\alpha)$ is $z$-unramified
and $\br{\theta_{D_1}(z)}\neq\br{\theta_{D_2}(z)}$,
and a {\it cold point} if each $\theta_{D_i}$ is $z$-ramified.
Then $\partial_z(\theta_{D_1})+\partial_z(\theta_{D_2})=0$ by \eqref{sumtozero}.

\Paragraph{\bf Splitting Ramification.}\label{splittingramification}
Each $\psi\in\H^1(F,\Z/n)$ determines a cyclic field extension, which we denote by $F(\psi)/F$.
In the situation of \eqref{setup} the normalization of $X$ in $L=F(\psi)$ is a regular relative
curve $Y/R$, and $Y\to X$ is flat, by \cite[Section 3]{BT11}.  Thus if $D\in\Div X$ is a prime divisor
then there are $g$ prime divisors $E_i\in\Div Y$ lying over $D$, each defining an extension $w_i$ of $v=v_D$
with ramification index $e=e(w_i/v)=|\partial_D(\psi)|$, 
and residue field $\k(E_i)=\k(D)((e\cdot\psi)(D))$ of degree $f$ over $\k(D)$,
such that $|\psi|=[F(\psi):F]=efg$.
We then have a 
commutative diagram
\[
\xymatrix{
&\H^2(F(\psi),\mu_n)\ar[r]^-{\partial_{E_i}}&\H^1(\k(E_i),\Z/n)\\
&\H^2(F,\mu_n)\ar[u]^\res\ar[r]_-{\partial_D}&\H^1(\k(D),\Z/n)\ar[u]_{e\cdot\res}\\
}\]
Conversely if $v_E$ is a valuation on $L$ determined
by a divisor $E\subset Y$ then since $Y\to X$ is flat,
$E$ lies over a divisor $D\subset X$, and $v_E$ extends $v_D$.

Continuing in the situation of \eqref{setup},
by definition an element $\beta\in\H^2(F(\psi),\mu_n)$ is unramified if $\partial_w(\beta)=0$ for all
discrete valuations $w$ on $F(\psi)$.
By purity for (regular) surfaces
it is enough to show $\partial_E(\beta)=0$ for each prime divisor $E\subset Y$, and then it
is in the image of the map $\H^2(Y,\mu_n)\to\H^2(F(\psi),\mu_n)$
(see \cite[(4.2)]{Br10}).  If $R=\Z_p$ then $\H^2(Y,\mu_n)=0$,
hence in this case if $\beta$ is unramified then it is zero (see \cite[Theorem 4.5]{Br10}).

\Paragraph{\bf Results from \cite{BMT13}.}
Assume the situation of \eqref{setup}.
Let $V=X-\D$.  This is not a scheme, but we use the notation heuristically and set
$\H^q(V):=\varinjlim\H^q(U)$, where the limit is over open subschemes
$U\subset X$ such that $X-U$ is supported in $\D$.
Since $V$ contains the generic points of $C$
the restriction map $\H^q(V)\longrightarrow\H^q(\k(C))$ is defined,
and as in \cite[(2.2)]{BMT13}, we have a commutative ladder
\begin{equation}\label{maplocseqs}
\xymatrix{
\cdots\ar[r]&\H^q(X)\ar[r]\ar[d]^\wr&\H^q(V)\ar[r]^-{\;\partial\;}\ar[d]&\bigoplus_{\D}\H^{q-1}(D,-1)
\ar[r]\ar[d]&\H^{q+1}(X)\ar[d]^\wr\ar[r]&\cdots\\
\cdots\ar[r]&\H^q(C)\ar[r]&\H^q(\k(C))\ar[r]^-{\;\delta\;}&
\bigoplus_{C_{(0)}}\H^{q+1}_z(C)\ar[r]&\H^{q+1}(C)\ar[r]&\cdots
}
\end{equation}
Let \[\Gamma_{\D_*}^q(\k(C))=\im(\H^q(V)\longrightarrow\H^q(\k(C)))\]
Then by \cite[Lemma 2.12]{BMT13}
$\Gamma_{\D_*}^q(\k(C))$ consists of tuples $\theta_C=(\theta_{C_i})$
where the $\theta_{C_i}\in\H^q(\k(C_i))$ {\it glue across $\S$ along $\D_*$},
that is, whenever $z\in\S$ is at the intersection of $C_i$ and $C_j$,
and $\pi_j$ and $\pi_i$ are the images in $\k(C_i)$ and $\k(C_j)$ of a local
equation $\pi_{D_z}\in\O_{X,z}$ for $D_z\in\D_*$, then the Witt decompositions of
$\theta_{C_i}$ and $\theta_{C_j}$ are
\begin{align}
\theta_{C_i,z}&=\theta_z^\circ+(\pi_j)\cdot\omega\;\in\H^q(K_{i,z})\label{glueatz}\\
\theta_{C_j,z}&=\theta_z^\circ+(\pi_i)\cdot\omega\;\in\H^q(K_{j,z})\notag
\end{align}
for some $\theta_z^\circ\in\H^q(\k(z))$ and $\omega\in\H^{q-1}(\k(z),-1)$.

Let $\H^q_\cs(V)=\ker(\H^q(V)\to\H^q(\k(C)))$,
let $\H^q_\cs(F)$ be its image in $\H^q(F)$,
and set $\H^q_\ns(F)=\H^q(F)/\H^q_\cs(F)$ (see \cite[Definition 2.11]{BMT13}).
The elements of $\H^q_\cs(F)$ are ``completely split'' in the sense that their images
in $\H^q(F_D)$ are zero for any prime divisor $D$ on $X$ (\cite[Lemma 2.12]{BMT13}),
and consequently the restriction maps $\H^q(F)\to\H^q(F_D)$, hence the residue and value maps, 
factor through $\H^q_\ns(F)$.

We require the following fundamental result for producing elements of $\H^q(F)$ with controlled
ramification.  Though our main application is the construction of elements of $\H^1(F,\Z/n)$, we state
the theorem in its entirety.

\Theorem{\cite[Theorem 2.15]{BMT13}.}\label{2.15}
Assume \eqref{setup}.
Then for all $q\geq 0$ there is a homomorphism
\[\lambda=\lambda_\D:\Gamma_{\D_*}^q(\k(C))\longrightarrow\H_\ns^q(F)\]
that fits into a commutative diagram
$$
\xymatrix{
\H^q_\ns(F)\ar[r]^-\partial&\bigoplus_{\D}\H^{q-1}(\k(D),-1)\\
\Gamma_{\D_*}^q(\k(C))\ar[r]^-\delta\ar[u]^{\lambda}&\bigoplus_{C_{(0)}}
\H^{q-1}(\k(z),-1)\ar[u]_\inf
}
$$
where $\delta$ is induced from \eqref{maplocseqs} and $\partial=\bigoplus_\D\partial_D$.
Let $\theta=\lambda(\theta_C)$, where $\theta_C=(\theta_{C_i})\in\Gamma_{\D_*}^q(\k(C))$ (so each $\theta_{C_i}$ is in $\H^q(\k(C_i))$). 
Then:
\begin{enumerate}[\rm a)]
\item
$\theta$ is defined on the generic points of $C$,
and $\theta(C_i)=\theta|_{\k(C_i)}=\theta_{C_i}$.
%$\res_{F|\prod_i F_i}\cdot\lambda =\inf_{\prod_i\k(C_i)|\prod_i F_i}:
%\Gamma_{\D_*}^q\to\bigoplus_i\H^q(F_i,\Lambda),$
\item
The ramification locus of $\theta$ is contained in $\D$.
\item
Suppose $D\in\D$ intersects $C$ at $z\in C_i$, $\pi_D\in\O_{X,z}$ is a local equation for $D$, 
$\bar\pi_D\in K_{i,z}$ is the image of $\pi_D$,
and $\theta_{C,z}=\theta^\circ+(\bar\pi_D)\cdot\omega$ is the corresponding Witt decomposition 
in $\H^q(K_{i,z})$
Then over $F_D$ we have the Witt decomposition
\[
\res_{F|F_D}(\theta)
=\inf_{\k(z)|\k(D)}(\theta^\circ)+(\pi_D)\cdot\inf_{\k(z)|\k(D)}(\omega)
\]
%Note that $\omega=\delta_z(\theta_C)$, maybe $\neq\partial_z(\theta_C)$ if $z\in\S$.
\item
If $\theta_C$ is unramified at a point $z\in C$,
then $\theta$ is unramified at any horizontal prime divisor $D$ lying over $z$, and 
$\theta(D)=\inf_{\k(z)|\k(D)}(\theta_C(z))$.
\end{enumerate}
\rm

\Remark
We will call any element $\theta\in\H^q(F)$ mapping to $\lambda(\theta_C)\in\H^q_\ns(F)$ a
{\it $\lambda$-lift of $\theta_C\in\Gamma_{\D_*}^q(\k(C))$}.

\section{Lemmas}

We begin with some preliminary results.
We will say that $k$ has no {\it special case} if $k(\mu_{2^{v_2(n)}})/k$ is cyclic.
This holds in particular if $n$ is odd or prime, or if $\car(k)>0$.
When $k$ has no special case Saltman's generalized Grunwald-Wang theorem \cite[Theorem 5.10]{Sa82}
produces for any $\k(C_i)$ and a finite set of local characters $\{\theta_{i,z}\in\H^1(K_{i,z},\Z/n):z\in C_i\}$ 
a global character $\theta_i\in\H^1(\k(C_i),\Z/n)$ of order $|\theta_i|=\lcm_z\{|\theta_{i,z}|]\}$
such that $\res_{\k(C_i)|K_{i,z}}(\theta_i)=\theta_{i,z}$.
We will call $\theta_i$ a (generalized) Grunwald-Wang lift of the $\theta_{i,z}$.

\Lemma\label{localstructure}
Assume the setup of \eqref{setup} such that $k$ has no special case.
Suppose $\alpha\in\H^2(F,\mu_n)$ 
ramifies at (only) two prime divisors $D_1$ and $D_2$ meeting transversely at a closed point $z\in X$.
Let $\pi_{D_1},\pi_{D_2}\in\O_{X,z}$ be local equations for $D_1$ and $D_2$ at $z$, and
suppose either {\rm (a)}
both $D_1$ and $D_2$ are vertical, or 
{\rm (b)}
each $\partial_{D_i}(\alpha)$ is $z$-unramified.
Then
\[\alpha=\alpha^\circ+(\pi_{D_1})\cdot\theta_1+(\pi_{D_2})\cdot\theta_2\]
for some $\alpha^\circ\in\H^2(\O_{X,z},\mu_n)$
and $\lambda$-lifts $\theta_i\in\H^1(F,\Z/n)$ with $D_{\theta_i}\subset\D$,
and in case {\rm (b)} $\theta_i\in\H^1(\O_{X,z},\Z/n)$.

\rm

\begin{proof}
We prove (a) first, setting $D_i=C_i$ for $i\in\{1,2\}$.
Set $\theta_{C_i}=\partial_{C_i}(\alpha)\in\H^1(\k(C_i),\Z/n)$.
Note that $\partial_z(\theta_{C_1})+\partial_z(\theta_{C_2})=0$ by \eqref{sumtozero},
and since $C_1$ and $C_2$ meet transversely at $z$
the image of $\pi_{C_i}$ in $\k(C_j)$ is a uniformizer for the valuation $v_z$ on $\k(C_j)$.
Let $D=\div(\pi_{C_1}-\pi_{C_2})\in\D_*$ be the distinguished divisor passing through $z$.
By the generalized Grunwald-Wang theorem of Saltman
(\cite[Theorem 5.10]{Sa82}) there exists an element
$\theta_{i,C}\in\Gamma^1_{\D_*}(\k(C))$ whose $C_i$-th component is $\theta_{C_i}$.
Let $\theta_i\in\H^1(F,\Z/n)$ be any lift of $\lambda(\theta_{i,C})\in\H^1_\ns(F,\Z/n)$ 
as in Theorem~\ref{2.15},
and set
\[\beta=(\pi_{C_1})\cdot\theta_1+(\pi_{C_2})\cdot\theta_2\]
Since the $\theta_i$ ramify only on $\D$ by Theorem~\ref{2.15}(b), 
$\beta$ is unramified at all divisors passing through $z$ different from $C_1$, $C_2$, and $D$.
Since $\theta_i$ is unramified at $C_i$ we compute $\partial_{C_i}(\beta)=\theta_i|_{\k(C_i)}$
using \eqref{cupresidue},
and this is $\theta_{C_i}\in\H^1(\k(C_i),\Z/n)$ by Theorem~\ref{2.15}(a).
Next, by \eqref{cupresidue} we compute
\[\partial_D(\beta)
=-(\bar\pi_{C_1})\cdot\partial_D(\theta_1)-(\bar\pi_{C_2})\cdot\partial_D(\theta_2)\]
where $\bar\pi_{C_i}$ is the image of $\pi_{C_i}$ in $\k(D)$.
We have $(\bar\pi_{C_1})=(\bar\pi_{C_2})$ by our choice of $D$,
and since $\partial_D(\theta_i)=\partial_z(\theta_{C_i})$ 
by Theorem~\ref{2.15}(c), $\partial_D(\theta_1)+\partial_D(\theta_2)=0$, hence $\partial_D(\beta)=0$.
We conclude $\alpha$ and $\beta$ have the same residues at divisors passing through $z$,
hence $\alpha^\circ\,\df\,\alpha-\beta\in\H^2(\O_{X,z},\mu_n)$ by injectivity and purity for surfaces 
(\cite[Theorem 7.2, Proposition 7.4]{AG60b}), as desired.

To prove (b) suppose first that $D_1=C_1$ and $D_2=C_2$ are both vertical.
Then (a) applies, and since each $\theta_{C_i}$ is $z$-unramified,
each $\theta_i$ is unramified at each divisor passing through $z$ by Theorem~\ref{2.15}, 
hence $\theta_i\in\H^1(\O_{X,z},\Z/n)$ by purity, as desired.

Suppose $D_i$ is horizontal for $i\in\{1,2\}$, and $z\in C_k$.
Then since $\theta_{D_i}$ is $z$-unramified and $\k(D_i)$ is complete with residue field $\k(z)$, 
$\theta_{D_i}$ is defined over $\k(z)$ by the Witt sequence.
Let $\psi_{C_k,z}=\theta_{D_i}\in\H^1(\k(z),\Z/n)\leq\H^1(K_{k,z},\Z/n)$,
let $\psi_{C_k}\in\H^1(\k(C_k),\Z/n)$ be any (generalized) Grunwald-Wang lift of $\psi_{C_k,z}$,
let $\psi_C\in\Gamma_{\D_*}^1(\k(C))$ be any element whose $C_k$-th component is $\psi_{C_k}$,
and let $\theta_i\in\H^1(F,\Z/n)$ be any lift of $\lambda(\psi_C)\in\H^1_\ns(F,\Z/n)$.
Then $\theta_i$ ramifies on $\D$ by Theorem~\ref{2.15}(b),
and since $\theta_{D_i}$ is $z$-unramified, $\partial_{D_i}(\theta_i)=0$ by Theorem~\ref{2.15}(c).
We compute $\partial_{D_i}((\pi_{D_i})\cdot\theta_i)=\theta_i|_{F_{D_i}}$, and
$\theta_i|_{F_{D_i}}=\theta_{D_i}$ by Theorem~\ref{2.15}(c).
Therefore $\partial_{D_i}((\pi_{D_i}\cdot\theta_i)=\theta_{D_i}$,
and if $D\neq D_i$ is any divisor running through $z$ then $\partial_D((\pi_{D_i})\cdot\theta_i)=0$.
Now if both $D_1$ and $D_2$ are horizontal then we define $\theta_1$ and $\theta_2$ as above.
If only $D_2$ is horizontal and $D_1=C_1$ is vertical 
then we choose $\theta_1=\lambda(\theta_{1,C})$ as in the proof of (1), 
where $\theta_{1,C}\in\Gamma_{\D_*}^1(\k(C))$ is any element whose $C_1$-th component is $\theta_{C_1}$.
In the latter case $\theta_1\in\H^1(\O_{X,z},\Z/n)$ since $\theta_{C_1}$ is $z$-unramified.
In either case
it follows immediately that $\alpha-(\pi_{D_1})\cdot\theta_1-(\pi_{D_2})\cdot\theta_2\in\H^2(\O_{X,z},\mu_n)$,
which proves (b).

Since the $\theta_i$ in all cases are $\lambda$-lifts we have $D_{\theta_i}\subset\D$ by
Theorem~\ref{2.15}(b), and in case (b) the $\theta_i$ are unramified through every divisor 
passing through $z$, hence $\theta_i\in\H^1(\O_{X,z},\Z/n)$ by purity and the Leray spectral sequence
(see \cite[Theorem 2.4]{BT10}).  This completes the proof.

\end{proof}

\Lemma\label{hotpoints}
Assume \eqref{setup}.
Suppose $n=\ell$ is prime, $\alpha\in\H^2(F,\mu_\ell)$, $D_\alpha\in\Div X$ has normal crossings,
and $z\in X$ is a hot point for $\alpha$, as in \eqref{hotpt}.
Then $\ell^2$ divides $\ind(\alpha)$.
\rm

\begin{proof}
Assume $z\in D_1\cap D_2$, with $D_i\subset D_\alpha$.
%Since $\ell$ is prime, by \cite[Proposition 7.2]{Br10} and \cite[Theorem 7.3]{Br10} 
By Lemma~\ref{localstructure}(b) we can write 
\[\alpha=\alpha^\circ+(\pi_{D_1})\cdot\theta_1+(\pi_{D_2})\cdot\theta_2\]
where $\alpha^\circ\in\H^2(\O_{X,z},\mu_\ell)$, 
$\theta_i\in\H^1(\O_{X,z},\Z/\ell)$ are $\lambda$-lifts, 
and $\pi_{D_i}\in\O_{X,z}$ is a local equation for $D_i$.
Set $\theta_{D_i}=\theta_i|_{F_{D_i}}$, then $\theta_{D_i}\in\H^1(\k(D_i),\Z/n)$ 
since $\partial_{D_i}(theta_{D_i})=0$ by the Witt sequence.
Since $z$ is a hot point $\theta_{D_i}$ is $z$-unramified, and we have values $\theta_{D_i}(z)$.
$\ell$ is prime we may assume $\theta_{D_2}(z)=0$ and $\theta_{D_1}(z)\neq 0$.

By the Nakayama-Witt index formula we have
\[\ind(\alpha_{F_{D_2}})=|\theta_{D_2}|\ind\left((\alpha^\circ+(\pi_{D_1})\cdot\theta_1)_{\k(D_2)(\theta_{D_2})}
\right)\]
Since $D_2\subset D_\alpha$ we have $\theta_{D_2}\neq 0$, hence $|\theta_{D_2}|=\ell$,
hence to show $\ell^2$ divides $\ind(\alpha)$
it suffices to prove that $(\alpha^\circ+(\pi_{D_1})\cdot\theta_1)_{\k(D_2)(\theta_{D_2})}$
is nontrivial.
Since $\theta_{D_2}(z)=0$, 
the valuation $v_z$ on $\k(D_2)$ determined by $z$ splits completely in $\k(D_2)(\theta_{D_2})$,
hence $\k(D_2)(\theta_2)$ has residue field $\k(z')=\k(z)$ with respect to any extension $v_{z'}|v_z$.
%Since $\partial_{D_2}(\alpha)=\theta_{D_2}$ is unramified and nontrivial, 
%$v_{D_2}$ extends uniquely to $F(\theta_2)$,
%hence $F_{D_2}(\theta_2)$ is a field extension with residue field $\k(D_2)(\theta_{D_2})$.
Since $D_\alpha$ has normal crossings at $z$, $v_{z'}(\bar\pi_{D_1})=1$, and we compute
using \eqref{cupresidue}
\[\partial_{z'}((\alpha^\circ+(\pi_{D_1})\cdot\theta_1)_{\k(D_2)(\theta_{D_2})})=
v_{z'}(\pi_{D_1})\cdot\theta_{D_1}(z)_{\k(z')}=\theta_{D_1}(z)\neq 0\]
Thus $\ell^2$ divides $\ind(\alpha_{F_{D_2}})$, hence $\ell^2$ divides $\ind(\alpha)$, as desired.
\end{proof}

The next lemma classifies distinguished divisors through a singular point $z\in\S$.

\Lemma\label{newDz}
Assume $C$ has normal crossings at the intersection $z\in C_i\cap C_j$,
and $\pi_i,\pi_j\in\O_{X,z}$ are local equations for $C_i$ and $C_j$.
Suppose a prime divisor $D$ on $X$ runs through $z$ with 
local equation $a_i\pi_i+a_j\pi_j\in\O_{X,z}$.
Then
\begin{enumerate}[\rm (a)]
\item
$D$ is regular at $z$ if and only if $(a_i,a_j)\O_{X,z}=\O_{X,z}$.
\item
$D$ is horizontal if and only if $a_i\not\in(\pi_j)$ and $a_j\not\in(\pi_i)$.
\item
$D$ intersects $C_i$ and $C_j$ transversely at $z$ if and only if $a_i,a_j\in\O_{X,z}^*$.
\end{enumerate}
\rm

\begin{proof}
We have the maximal ideal $\frak m_z=(\pi_i,\pi_j)\O_{X,z}$ since $C$ has normal crossings at $z$.
Suppose given $a_{11},a_{21}\in\O_{X,z}$ such that $(a_{11},a_{21})\O_{X,z}=\O_{X,z}$.
Then there exist $a_{12},a_{22}\in\O_{X,z}$ such that $a_{11}a_{22}-a_{21}a_{12}=1$, so 
the matrix
\[A=\begin{pmatrix}a_{11} & a_{12}\\ a_{21} & a_{22}\end{pmatrix}\]
is invertible.
It follows that 
$\frak m_z=(\pi_i',\pi_j')\O_{X,z}$ for $\pi_i'=a_{11}\pi_i+a_{21}\pi_j$ and $\pi_j'=a_{12}\pi_i+a_{22}\pi_j$
by the invertibility of $A$, and
$\pi_i'$ is regular as part of a regular system
of generators for $\frak m_z$. 
Thus $D=\div(\pi_i')$ is regular at $z$ by definition.
Conversely if $D$ is regular at $z$ then locally $D\cap\Spec\O_{X,z}=\div(\pi_i')$
for a regular element $\pi_i'=a_{11}\pi_i+a_{21}\pi_j$, and if $\pi_j'=a_{12}\pi_i+a_{22}\pi_j$
completes the regular system at $z$ then we obtain an invertible matrix $A$ as above,
and the condition $\det(A)\in\O_{X,z}^*$ shows $(a_{11},a_{21})\O_{X,z}=\O_{X,z}$.
This proves (a).

Since $D$ is a prime divisor it is either horizontal or vertical, 
and since $C$ has normal crossings at $z$, 
$D$ is horizontal if and only if does not coincide with $C_i$ or $C_j$ at $z$,
i.e., $a_i\pi_i+a_j\pi_j\not\in(\pi_i)\cup(\pi_j)$.
Equivalently $a_i\pi_i\not\in(\pi_j)$ and $a_j\pi_j\not\in(\pi_i)$,
i.e., $a_i\not\in(\pi_j)$ and $a_j\not\in(\pi_i)$.
This proves (b).

Finally,
$D$ is transverse to both $C_i$ and $C_j$ if and only if
$(a_i\pi_i+a_j\pi_j,\pi_i)\O_{X,z}=(a_j\pi_j,\pi_i)\O_{X,z}=\frak m_z$ and $(a_i\pi_i+a_j\pi_j,\pi_j)\O_{X,z}=(a_i\pi_i,\pi_j)\O_{X,z}=\frak m_z$,
which is equivalent to $a_i,a_j\in\O_{X,z}^*$.
\end{proof}

\Remark
By Lemma~\ref{newDz}(c) we may choose for $D_z\in\D_*$ the divisor associated 
to any $\O_{X,z}^*$-linear combination of $\pi_i$ and $\pi_j$.

We use the next lemma to glue across cold points.

\Lemma\label{step1}
Assume the setup of \eqref{setup},
$n=\ell$ is prime, $\alpha\in\H^2(F,\mu_\ell)$, $D_\alpha$ has normal crossings on $X/R$,
and $z$ is a cold point for $\alpha$ at the intersection of vertical components $C_1,C_2\subset D_\alpha$.
Then there exists a regular horizontal divisor $D$ running through $z$ and
transverse to $C_1$ and $C_2$, such that $\partial_{C_1}(\alpha)$ and $-\partial_{C_2}(\alpha)$
glue at $z$ along $D$ as in \eqref{glueatz}.
\rm

\begin{proof}
Since $D_\alpha$ has normal crossings, $C_1$ and $C_2$ meet transversely at $z$, and
we have maximal ideal $\frak m_z=(\pi_1,\pi_2)\O_{X,z}$ where $\pi_i$ is a local equation for $C_i$ at $z$.
By Lemma~\ref{newDz} the local equation $\pi_{D_z}:=\pi_1+\pi_2\in\O_{X,z}$ defines a distinguished
divisor $D_z\in\D_*$.
Set $\theta_{C_i}:=\partial_{C_i}(\alpha)$.
The Witt decomposition of $\theta_{C_i}$ at $z$ along $D$ as in \eqref{glueatz} is
$\theta_{C_i}^\circ+(\bar\pi_j)\cdot\partial_z(\theta_{C_i})\in\H^1(K_{i,z},\Z/\ell)$
for some $\theta_{C_i}^\circ\in\H^1(\k(z),\Z/\ell)$, where $\bar\pi_j$ is the image of $\pi_j$ in $K_{i,z}$.
Set $\omega:=\partial_z(\theta_{C_1})$. 
Then $\omega$ has order $\ell$ since $z$ is a cold point and $\ell$ is prime,
and $\partial_z(\theta_{C_2})=-\omega$
since $\partial_z(\theta_{C_1})+\partial_z(\theta_{C_2})=0$ (by \eqref{sumtozero}).
Thus $\mu_\ell\subset\k(z)$, hence $\theta_{C_i}^\circ$ is a Kummer character, of the form 
$(\bar a_i)\cdot\omega$ for some $\bar a_i\in\k(z)^*$.
Choose preimages $a_i\in\O_{X,z}^*$, and let $b_1=a_1^{-1}$.
Then the divisor $D'=\div(b_1\pi_1+a_2\pi_2)$ is transverse to both $C_1$ and $C_2$ by Lemma~\ref{newDz},
and the Witt decompositions of the $\theta_{C_i}$ at $z$ along $D'$ are
\begin{align*}
\theta_{C_1,z}&=\theta_{C_1}^\circ-(\bar a_2)\cdot\omega+(\bar a_2\bar\pi_2)\cdot\omega
=(\bar a_1\bar a_2^{-1})\cdot\omega+(\bar a_2\bar\pi_2)\cdot\omega
\in\H^1(K_{1,z},\Z/\ell)\\
\theta_{C_2,z}&=\theta_{C_2}^\circ+(\bar b_1)\cdot\omega-(\bar b_1\bar\pi_1)\cdot\omega
=(\bar a_2\bar a_1^{-1})\cdot\omega-(\bar b_1\bar\pi_1)\cdot\omega
\in\H^1(K_{2,z},\Z/\ell)\\
\end{align*}
where $\bar\pi_j$, $\bar b_1$, and $\bar a_j$ are the images in $K_{i,z}$.
Thus $\theta_{C_1}$ and $-\theta_{C_2}$ glue at $z$ along $D'$ by \eqref{glueatz}.
\end{proof}

We next show how to break cycles in the dual graph $\G_{D_\alpha}$ by blowing up.
A {\it chordless cycle} (or {\it hole}) of a graph $\G$ is a sequence of vertices of $\G$
such that each pair of adjacent vertices are connected by an edge in $\G$, and no non-adjacent
vertices are connected by an edge.

\Lemma\label{step2}
Assume the setup of \eqref{setup},
$n=\ell$ is prime, $\alpha\in\H^2(F,\mu_\ell)$, $D_\alpha$ has normal crossings on $X/R$,
$z$ is at the intersection of vertical components $C_1,C_2\subset D_\alpha$,
and $z$ is neither a hot point or a cold point for $\alpha$.
Then there exists a blowup $X'$ of $X$ centered at $z$ such that 
$\alpha$ is unramified on some 
irreducible component $E$ of the exceptional fiber.
In particular if $\b_1(D_\alpha)\geq 1$ and $z\in D_\alpha$ is a closed point in $\S$
that corresponds in $\G_{D_\alpha}$ to an edge of a chordless cycle,
and $z$ is not a hot or cold point for $\alpha$,
then there exists a blowup $X'$ of $X$ over which the divisor $D'_\alpha\in\Div X'$ has Betti
number $\b_1(D'_\alpha)<\b_1(D_\alpha)$.
\rm

\begin{proof}
Set $\theta_{C_i}=\partial_{C_i}(\alpha)$.
Since $z$ is not a cold point each $\theta_{C_i}$ has a value $\theta_{C_i}(z)$,
and since $z$ is not hot we have $\theta_{C_2}(z)=n_z\theta_{C_1}(z)$ for some $n_z\in(\Z/\ell)^*$.
Let $E$ be the exceptional divisor of the blowup of $X$ at $z$.
By Lemma~\ref{localstructure}(b) we may write 
$\alpha=\alpha^\circ+(\pi_{C_1})\cdot\theta_1+(\pi_{C_2})\cdot\theta_2$
with $\theta_i\in\H^1(\O_{X,z},\Z/\ell)$.
Since $v_E(\pi_{D_i})=1$ for each $i$, we compute
$\theta_E=\partial_E(\alpha)=(\theta_1+\theta_2)|_{F_E}$ using \eqref{cupresidue}, 
and since $\partial_E$ factors through $\k(z)$ this
is $\theta_{C_1}(z)+\theta_{C_2}(z)=(n_z+1)\theta_{C_1}(z)$.
Thus after finitely many blowups of $X$ we reach an exceptional divisor $E'$ at which $\theta_{E'}=0$,
proving the first statement.

Assume now that moreover $z\in\S$ corresponds in $\G_{D_\alpha}$ to an edge of a chordless cycle.
Let $X'\to X$ be the composite blowup (centered at $z$), and let $D_\alpha'$ be 
the divisor of $\alpha$ on $X'$.
Clearly $D_\alpha'$ has normal crossings on $X'$.
The effect of a blowup on a dual graph in general is to divide an edge and its two vertices into
two edges and three vertices, hence blowups preserves cycles on divisors. 
However since $\theta_{E'}=0$ we have removed
a vertex from the blowup of $\G_{D_\alpha}$, thus breaking the chordless cycle to which
$z$ belonged, without creating any new cycles.  Thus $\b_1(D_\alpha')<\b_1(D_\alpha)$.
\end{proof}

\section{Main Theorem}

Assume the setup of \eqref{setup}, $n=\ell$ is prime, $\alpha\in\H^2(F,\mu_\ell)={}_\ell\Br(F)$,
and $D_\alpha$ has normal crossings on $X/R$.
We may assume the horizontal components of $D_\alpha$ are in $\D_\S$, by blowing up $X$ if necessary.
Let $C_\alpha=D_\alpha-\D(D_\alpha)$, which is the ``vertical'' part of $D_\alpha$.
Then $C_\alpha$ has normal crossings, hence has a dual graph $\G(C_\alpha)$,
and since all horizontal components of $D_\alpha$ are in $\D_\S$ it is clear that 
$\b_1(D_\alpha)=\b_1(C_\alpha)$.
We introduce some terminology.
Recall a {\it chordless cycle} (or {\it hole}) of a graph $\G$ is a sequence of vertices of $\G$
such that each pair of adjacent vertices are connected by an edge in $\G$, and no non-adjacent
vertices are connected by an edge.
\begin{itemize}
\item
A {\it tree} in $C_\alpha$ is a connected subset whose image in $\G(C_\alpha)$ is a tree.
\item
An {\it isolated tree} in $C_\alpha$ is a maximal connected component of $C_\alpha$ that is a tree.
\item
An {\it isolated-tree point} is a singular point on an isolated tree.
\item
A {\it cycle} in $C_\alpha$ is a subset whose image in $\G(C_\alpha)$ is a chordless cycle.
\item
A {\it cycle point} is a singular point of $C_\alpha$ whose corresponding edge 
in $\G(C_\alpha)$ is part of a cycle.
\item
A {\it cycle cluster} in $C_\alpha$ is a set of cycles,
maximal with respect to the property
that one may travel from one cycle in the cluster to another on the components of cycles.
\item
A {\it connecting path} in $C_\alpha$ is a maximal connected set of irreducible components that are not components of
cycles, and the set intersects more than one cycle cluster.
\item
A {\it connecting point} of $C_\alpha$ is a singular point of a connecting path.
\item
A {\it tail} in $C_\alpha$ is a maximal connected set of irreducible components that are not components of
cycles, and the set intersects exactly one cycle cluster at a single point.
\item
A {\it tail point} of $C_\alpha$ is a singular point of a tail.
\end{itemize}
With this terminology, $C_\alpha$ is a union of isolated trees, cycle clusters, connecting paths, and tails, and
the singular points of $C_\alpha$
are either isolated-tree points, cycle points, connecting points, or tail points.

\Theorem\label{cyclic}
Assume \eqref{setup},
$\alpha\in\Br(F)$ has prime period $\ell\neq p$.
Then there exists a $\Z/\ell$-cyclic extension $L/F$ such that $\alpha_L$ is unramified.
\rm

\begin{proof}
Let $X$ be a regular model for $F$ as in \eqref{setup}.
We will use the notation $\theta_D:=\partial_D(\alpha)\in\H^1(\k(D),\Z/\ell)$ for $D\in\Div X$,
and if $D=C_i$ for some irreducible component $C_i$ of $C$ we denote by
$\theta_{C_i,z}$ the image of $\theta_{C_i}$ in $\H^1(K_{i,z},\Z/\ell)$, 
where $K_{i,z}=\Frac(\O_{C_i,z}^\h)$ is as in \eqref{setup}.

By blowing up $X$ if necessary
we may assume that $D_\alpha\cup C$ has normal crossings, 
each singular point of $\Supp(D_\alpha)$
lies on a vertical component of $\Supp(D_\alpha)$,
all horizontal components of $D_\alpha$ are in $\D_\S$,
and the dual graph $\G_C$ is bipartite, so that 
$C$ is union of two disjoint sets of irreducible components $C^+$ and $C^-$.
We call any such model $X$ {\it $\alpha$-acceptable}.
Since the blowup at a closed point preserves normal crossings of divisors,
an even number of blowups of an $\alpha$-acceptable model is again $\alpha$-acceptable.

For every $\alpha\in\H^2(F,\mu_\ell)$ there exists an $\alpha$-acceptable model $X/R$ over which
$\b_1(D_\alpha)$ is minimal, among $\alpha$-acceptable models.
To prove the theorem we will induct on
the minimum value assumed by $\b_1(D_\alpha)$ on any $\alpha$-acceptable model $X/R$.  
Recall $\b_1(D_\alpha)$ is the number of chordless cycles in $\G_{D_\alpha}$, which is the
same as the number of loops in $D_\alpha$.
Since a blowup cannot join two disconnected components it does not affect $\b_1(D_\alpha)$.
We are therefore free to blow up a model $X$ over which $\G_C$ is defined
and $\b_1(D_\alpha)$ is minimal until we obtain a model that is additionally 
$\alpha$-acceptable.

\Paragraph\label{treecase}
{\it Inductive procedure for trees.}
Suppose $\alpha\in\H^2(F,\mu_\ell)$ and there exists an $\alpha$-acceptable
model $X$ such that $\b_1(D_\alpha)=0$, i.e., $D_\alpha$ is a tree.
Let $C_\alpha$ be the vertical components of $D_\alpha$ as above,
then $\b_1(C_\alpha)=0$.
If $C_\alpha\neq\varnothing$ then sequence the components of a connected component
of $C_\alpha$ by $C_1,\dots,C_r$, so that $C_d\cap(\bigcup_{i<d} C_i)\neq\varnothing$.
Fix such a connected component,
and define $\psi_{C_1}$ to be any prime-to-$\ell$ multiple of $\theta_{C_1}$,
and if $r\neq 1$ then inductively define $\psi_{C_d}$ for $d:1<d\leq r$ as follows.
\begin{enumerate}[(I)]
\item
If $i<d$, $z\in C_i\cap C_d$, and $\theta_{C_i,z}$ is $z$-ramified then
let $\psi_{C_d}=\theta_{C_d}$ if $C_d\subset C^+$ and $\psi_{C_d}=-\theta_{C_d}$ if $C_d\subset C^-$.
Then $\partial_z(\psi_{C_i})=\partial_z(\psi_{C_d})$.
\item
If $i<d$, $z\in C_i\cap C_d$, and $\theta_{C_i,z}$ is $z$-unramified, 
let $\psi_{C_d}=n_z\theta_{C_d}$ where $n_z$ is a prime-to-$\ell$ 
number such that $\theta_{C_i}(z)=n_z\theta_{C_d}(z)$.
The number $n_z$ exists since $D_\alpha$ has no hot points.
%since $\b_1(C_\alpha)=0$, $\psi_{C_d}$ for $d>i$ will not have been determined by
%an earlier step.
\end{enumerate}
Assign $\psi_{C_i}$ in this way for each connected component of $C_\alpha$.
Suppose $C_\alpha$ has $s$ irreducible components, and let
$C_{s+1},\dots,C_m$ denote some ordering of the remaining irreducible components of $C$.
Note if $i\leq s<d$, and $z\in C_d\cap C_i$ then $\psi_{C_i}$ is $z$-unramified
since $\theta_{C_d}=0$ and $\partial_z(\theta_{C_i}+\theta_{C_d})=0$ (by \eqref{sumtozero}).
Now inductively define (local) data $\psi_{C_d,z}\in\H^1(K_{d,z},\Z/\ell)$ for 
$d>s$ and various $z\in C_d$ as follows.
\begin{enumerate}[(A)]
\item
If $d>i$ and $z=C_d\cap C_i$ set $\psi_{C_d,z}=\psi_{C_i,z}=\psi_{C_i}(z)\in\H^1(\k(z),\Z/\ell)$.
\item
If $d<i$ and $z=C_d\cap C_i$ set $\psi_{C_d,z}=0$.
\item
If $z\in C_d\cap\D(D_\alpha)$ 
let $\psi_{C_d,z}=\theta_D\in\H^1(\k(z),\Z/\ell)\leq\H^1(K_{d,z},\Z/\ell)$,
where $D\in\D(D_\alpha)$ passes through $z$.
\end{enumerate}
Note in (C) that indeed $\theta_D\in\H^1(D,\Z/\ell)\isom\H^1(\k(z),\Z/\ell)$, and that
$D\in\D_\S$ since we have assumed $\D(D_\alpha)=\D_\S(D_\alpha)$.
Now let $\psi_{C_d}\in\H^1(\k(C_d),\Z/\ell)$ be any element with images the $\psi_{C_d,z}$
from $(A,B,C)$.  
Such an element exists by Saltman's generalized Grunwald-Wang Theorem \cite[Theorem 5.10]{Sa82}, 
since here there is no special case and there are finitely many singular points $z$ on $C\cup D_\alpha$.

Finally, if $C_i$ and $C_j$ are in $C_\alpha$ and intersect at $z\in\S$
and $\theta_{C_i}$ (hence $\theta_{C_j}$) is $z$-ramified,
then by Lemma~\ref{step1} and (I) there exists a choice $D=D_z\in\D_*$ such
that $\psi_{C_i}$ and $\psi_{C_j}$ glue at $z$ along $D$.
Then by (II) the $\psi_{C_k}$ glue along all singular points of $C_\alpha$,
and by (A,B) the $\psi_{C_k}$ glue along all other points of $\S$.
By \eqref{glueatz} there exists
an element $\psi_C\in\Gamma^1_{\D_*}(\k(C))\leq\H^1(\k(C),\Z/\ell)$, hence a $\lambda$-lift
$\psi\in\H^1(F,\Z/\ell)$ of $\psi_C$ by Theorem~\ref{2.15}.

\Paragraph\label{splitting}
{\it Splitting ramification for trees.}
Let $L=F(\psi)$ and let $Y$ be the normalization of $X$ in $L$.  
We claim that $\alpha_L$ is unramified.
By \eqref{splittingramification} it is enough to show that for each $D\subset D_\alpha$,
and $E\in\Div Y$ lying over $D$, $e\cdot\res_{\k(D)|\k(E)}(\theta_D)=0$,
where $e=|\partial_D(\psi)|$ and $\k(E)=\k(D)((e\cdot\psi)(D))$.
Suppose then that $D\subset D_\alpha$, and $E\in\Div Y$ lies over $D$.  
Note that $D\in C\cup\D_\S$.  Since $\psi$ is a $\lambda$-lift we may use Theorem~\ref{2.15}
to analyze its ramification behavior.
\begin{enumerate}[{\rm (a)}]
\item
If $D=C_i\subset C$ then $\partial_D(\psi)=\partial_{C_i}(\psi)=0$ by Theorem~\ref{2.15}(b),
and we have $\br{\psi_{C_i}}=\br{\theta_{C_i}}$ by (I) and (II).
Therefore $e=|\partial_{C_i}(\psi)|=1$ and $\k(E)=\k(C_i)(\psi(C_i))$,
hence $\k(E)=\k(C_i)(\psi_{C_i})=\k(C_i)(\theta_{C_i})$ by Theorem~\ref{2.15}(a), 
hence $\res_{\k(C_i)|\k(E)}(\theta_{C_i})=0$, 
hence $\partial_E(\alpha_L)=0$.
\item
If $D\in\D_\S$ and $\partial_D(\psi)\neq 0$ then $e=|\partial_D(\psi)|=\ell$ since $\ell$ is prime,
and since $\res_{\k(D)|\k(E)}(\theta_D)$ has order dividing $\ell$, $\partial_E(\alpha_L)=0$.
\item
If $D\in\D_\S$, $\partial_D(\psi)=0$, $D\cap C=z$ is on $C_i$, and $\theta_{C_i}=0$,
then $e=1$ so $\k(E)=\k(D)(\psi(D))$, and this is $\k(D)(\psi_{C_i}(z))$ by Theorem~\ref{2.15}(d),
where we identify $\psi_{C_i}(z)\in\H^1(\k(z),\Z/\ell)$ with its image in $\H^1(\k(D),\Z/\ell)$.
Since $\partial_D(\psi)=0$ we have $\partial_z(\psi_{C_i})=0$ by Theorem~\ref{2.15}(c),
hence $\psi_{C_i}(z)=\psi_{C_i,z}$,
and we have $\psi_{C_i,z}=\theta_D\in\H^1(\k(z),\Z/\ell)$ by (C), since $\theta_{C_i}=0$. 
Therefore $\k(E)=\k(D)(\theta_D)$, and $\partial_E(\alpha_L)=0$.
\item
If $D\in\D_\S$, $\partial_D(\psi)=0$, $D\cap C=z$ is on $C_i$,
and $\theta_{C_i}\neq 0$, then $e=1$ so $\k(E)=\k(D)(\psi(D))$, 
and $\psi(D)=\psi_{C_i}(z)$ by Theorem~\ref{2.15}(d) 
(where again we identify $\psi_{C_i}(z)\in\H^1(\k(z),\Z/\ell)$ with its image in $\H^1(\k(D),\Z/\ell)$).
%so $\k(E)=\k(D)(\psi_{C_i,z})$.

Since $\theta_{C_i}\neq 0$ we have $\br{\psi_{C_i}}=\br{\theta_{C_i}}$ by (I) and (II),
and since $\partial_D(\psi)=0$ we have 
$\partial_z(\psi_{C_i})=0$ by Theorem~\ref{2.15}(c), hence $\partial_z(\theta_{C_i})=0$.
Thus both $\theta_D$ and $\theta_{C_i}$ are nonzero and $z$-unramified.  
Since $D_\alpha$ has no hot points, 
$\br{\theta_{C_i}(z)}=\br{\theta_D(z)}$, and $\theta_D(z)=\theta_D$ since $\k(D)$ is already complete.
Therefore $\br{\psi_{C_i}(z)}=\br{\theta_D}\leq\H^1(\k(z),\Z/\ell)$.
We conclude $\k(E)=\k(D)(\theta_D)$,
hence $\res_{\k(D)|\k(E)}(\theta_D)=0$, and once more $\partial_E(\alpha_L)=0$.
\end{enumerate}
This completes the proof that when $\b_1(D_\alpha)=0$ for an $\alpha$-acceptable model $X$, 
there exists a $\Z/\ell$-cyclic extension $L/F$ such that $\alpha_L$ is unramified.

\Paragraph\label{generalcase}
{\it Case $\b_1(D_\alpha)>0$.}
Assume we have shown that for any $\alpha'\in\H^2(F,\Z/\ell)$ for which there is an
$\alpha'$-acceptable model $X'$ such that $\b_1(D_{\alpha'})\leq N$, there exists
a cyclic extension $L'/F$ such that $\alpha'_{L'}$ is unramified.
Fix $\alpha\in\H^2(F,\mu_\ell)$, and suppose that the minimum value of $\b_1(D_\alpha)$
over all $\alpha$-acceptable models is $\b_1(D_\alpha)=N+1$.
We must show there exists a $\Z/\ell$-cyclic extension $L/F$ such that $\alpha_L$ is
unramified.

Let $X$ be an $\alpha$-acceptable model over which $\b_1(D_\alpha)$ is minimal.
(Note as before we have $\b_1(D_\alpha)=\b_1(C_\alpha)$ since $\D(D_\alpha)\subset\D_\S$.)
Then every cycle point for $D_\alpha$ is cold by the second statement of Lemma~\ref{step2}, 
since otherwise by blowing up we could reduce $\b_1(D_\alpha)$ on another $\alpha$-acceptable
model.
Similarly we may assume that every connecting point is cold,
since otherwise we may blow up until the connecting path
becomes two tails by the first statement of Lemma~\ref{step2}.
By Lemma~\ref{step1} we may choose $D_z\in\D_*$ for every cold 
cycle point or connecting point $z\in\S$ on $C_\alpha$
so that if $z\in C_i\cap C_j$, then  
$\theta_{C_i}$ and $-\theta_{C_j}$ glue at $z$ along $D_z$ as in \eqref{glueatz}. 
Since $C_\alpha$ is bipartite we may set $\psi_{C_i}=\pm\theta_{C_i}$ accordingly, 
as in (I).
Thus the $\psi_{C_i}$ for the irreducible components $C_i\subset C_\alpha$
of cycle clusters or connecting paths glue at their singular points along $\D_*$.
The remaining components of $C_\alpha$ are isolated trees or tails, 
and for them we may define the $\psi_{C_i}$ inductively as we did for trees in \eqref{treecase}(I,II).
Thus we obtain elements $\psi_{C_i}$ for every $C_i\subset C_\alpha$, that
glue over the singular points of $C_\alpha$ along $\D_*$.
Finally, we extend to the remaining components of $C$
by defining the $\psi_{C_i}$ as in \eqref{treecase}(A,B,C),
so that the $\psi_{C_i}$ are compatible with the $\theta_D$ for isolated
components $D\subset D_\alpha$, and all components glue across singular points with those 
$\psi_{C_j}$ already defined.

Since the $\psi_i$ glue along all $z\in\S$ along $\D_*$ we have an element $\psi_C\in\Gamma^1_{\D_*}(\k(C))$,
and we let $\psi\in\H^1(F,\Z/\ell)$ be a $\lambda$-lift of $\psi_C$, and set
$L=F(\psi)$.
We claim that $\alpha_L$ is unramified.
As in the case \eqref{splitting} for trees, 
by \eqref{splittingramification} it is enough to show that for each $D\subset D_\alpha$
and $E\subset Y$ lying over $D$, $e\cdot\res_{\k(D)|\k(E)}(\theta_D)=0$,
where $e=|\partial_D(\psi)|$ and $\k(E)=\k(D)((e\cdot\psi)(D))$.
As before we have $D\subset C\cup\D_\S$.
\begin{enumerate}[{\rm (i)}]
\item 
If $D=C_i$ and $C_i$ is part of a tail or an isolated tree, 
or if $D\in\D_\S$ and $D$ intersects $C$ at such a component,
then the computation $\partial_E(\alpha_L)=0$ follows from \eqref{splitting}(a,b,c,d).
\item
If $D=C_i$ is part of a cycle cluster or connecting path then
we have $\partial_E(\alpha_L)=0$ by \eqref{splitting}(a).
\item
If $D\in\D_\S$ and $\partial_D(\psi)\neq 0$ the $\partial_E(\alpha_L)=0$ by \eqref{splitting}(b).
\item
If $D\in\D_\S$, $\partial_D(\psi)=0$, and $D\cap C=z$ lies on a component $C_i$
at which $\theta_{C_i}=0$, then $\psi_{C_i}$
is defined as in \eqref{treecase}(C), and $\partial_E(\alpha_L)$ follows from \eqref{splitting}(c).
\item
If $D\in\D_\S$, $\partial_D(\psi)=0$, and $D\cap C=z$ lies on a component $C_i$
of a cycle cluster or connecting path, then $\theta_{C_i}\neq 0$, and $\psi_{C_i}=\pm\theta_{C_i}$.
Therefore we are in the same situation as \eqref{splitting}(d), and we conclude that
$\partial_E(\alpha_L)=0$.
\end{enumerate}
We conclude that $\alpha_L$ is unramified.
The result now follows by induction.

\end{proof}

\Remark
This result compares with \cite[Theorem 7.13]{Sa08}, which applies to a general regular surface
but only when $F$ contains a primitive $\ell$-th root of unity.

We now reprove \cite[Theorem 5.1, Corollary 5.2]{Sa07}

\Corollary\label{cyclicQp}
{\rm (Cyclicity in Prime Degree)}
If $F$ is the function field of a smooth $p$-adic curve and $A$ is an $F$-division
algebra of index $\ell\neq p$, then $A$ is cyclic.
\rm

\begin{proof}
This is immediate since the $L$ of Theorem~\ref{cyclic} has trivial unramified Brauer group
(see e.g. \cite[Theorem 4.5]{Br10}),
hence $L$ splits the class $[A]\in\H^2(F,\mu_\ell)$.
\end{proof}

\Corollary\label{hotpointcriterion}
{\rm (Hot Point Criterion)}
If $F$ is the function field of a smooth $p$-adic curve,
$\alpha\in\Br(F)$ has prime period $\ell\neq p$, and $D_\alpha$ has normal
crossings, then $\ind(\alpha)=\ell$ if and only if $\alpha$ has no hot points.
\rm

\begin{proof}
If $\alpha$ has a hot point then $\ind(\alpha)\neq\ell$ by Lemma~\ref{hotpoints}.
Conversely if $\alpha$ has no hot points then $\ind(\alpha)=\ell$ by 
Corollary~\ref{cyclicQp}.

\end{proof}

\Remark
In the situation of Corollary~\ref{hotpointcriterion} we have $\ind(\alpha)|\ell^2$ by \cite{Sa97}
(see also \cite[Corollary 5.2]{BMT13}), hence $\ind(\alpha)=\ell^2$ if and only if $\alpha$ has a hot point.

\bibliographystyle{abbrv} %other choices are plain or abbrv or alpha
\bibliography{hnx.bib}

\end{document}